\documentclass{amsart}

\usepackage{verbatim,amsmath,latexsym,amssymb,amsbsy}

\input xy
\xyoption{all}

\newtheorem{prop}{Proposition}[section]
\newtheorem{lem}[prop]{Lemma}
\newtheorem{carol}[prop]{Corollary}
\newtheorem{theo}[prop]{Theorem}

\theoremstyle{remark}
\newtheorem{hyp}[prop]{Hypothesis}

\newtheorem{example}[prop]{Example}

\def\AA{{\mathbb A}}     
\def\CC{{\mathbb C}}

\def\PP{{\mathbb P}}     
   
\def\RR{{\mathbb R}} 
\def\ZZ{{\mathbb Z}}

\def\calV{{\mathcal V}}
\def\calW{{\mathcal W}}

\def\Mtwo{{\em Macaulay} 2\expandafter}

\numberwithin{equation}{section}
\begin{document}
\title{A result about Picard-Lefschetz Monodromy}
\author{ Darren Salven Tapp	}
\address{Department of Mathematics,
Purdue University,
150 North University Street,
West Lafayette, IN  47907-2067}

\begin{abstract}
Let $f$ and $g$ be reduced homogeneous polynomials in separate sets of variables.  We establish 
a simple formula that relates the eigenspace decomposition of the monodromy operator 
on the Milnor fiber cohomology of $fg$ to that of $f$ and $g$ separately.    
We use a relation between local systems and Milnor fiber cohomology that has been established by D. Cohen and A. Suciu.
\end{abstract}

\maketitle

\section{Brief Introduction and Statement of Results}

Thom-Sebastiani type Theorems have a rich history.  This study was initiated by
Sebastiani and Thom \cite{article:Thom} and improved by
others \cite{article:Oka,incollection:Sakamoto}.  %There may be more 
Perhaps the most
successful generalization has been achieved by N\'emethi
\cite{article:Nemethi:gen}.  He considers the germs of three 
holomorphic functions $f,g,p$ at the origin of $\CC^n, \CC^m$ and $\CC^2$
respectively, and then draws conclusions about the topology of the
Milnor fiber \cite{book:Milnor} of $p(f,g)$.   N\'emethi discovered an expression for the Weil zeta function 
of $p(f,g)$ in terms of the monodromy representations of $f$ and $g$ as well as the several variable Alexander polynomial of $p$ 
\cite{article:Nemethi:zeta,article:Nemethi:genzeta}.

In this paper we will investigate the case when $f$ and $g$
are homogeneous polynomials and $p(x,y) = xy$.  We will
construct a fibration different from those of 
\cite[Theorem 2]{incollection:Sakamoto}, and \cite{article:Nemethi:gen,article:Nemethi:zeta,article:Oka}.  
We will then produce a formula for the
eigenspace decomposition of the Picard-Lefschetz monodromy of $fg$ in
terms of those of $f$ and $g$. 

The following describes the situation we consider,
\begin{hyp}\label{hyp:assumptions}
Let $f \in \CC[x_1, \ldots , x_n]$ and $ g \in \CC[y_1,\ldots, y_m]$ 
be homogeneous and reduced of positive degrees $r$ and $s$ respectively.
\end{hyp}
In this case the restriction $f : \CC^n\setminus \{ x \mid f(x) = 0 \} \to \CC^*$ defines a fibration.  The fiber of this 
fibration is called the Milnor fiber of $f$, denoted $F_f = f^{-1}(1)$ .  
Lifting the path $t \mapsto \exp (2\pi i t): 0\leq t \leq 1$ 
in $\CC^*$ induces, through a local trivialization of $f$, a diffeomorphism of the fiber.  This map 
will be called a \emph{geometric} Picard-Lefschetz(PL) monodromy of $f$.  
A geometric PL monodromy of $f$ 
induces a map on the cohomology algebra $H^*(F_f,\CC)$ which we will call the \emph{algebraic} PL monodromy of $f$.  
In a more general setting we may have any smooth fibration $F\to M \to N$.  
Given a loop in the base space $N$ we may again 
construct a diffeomorphism of the fiber; in the case when this diffeomorphism is homotopic to the identity for any loop we choose, 
we say that the fibration has \emph{trivial} geometric monodromy.

For any $M \subseteq \CC^t
\setminus \{0\}$ we denote by $M^*$ the image of $M$ under the Hopf
fibration 
\[
\rho:\CC^t \setminus \{0\} \to \PP^{t-1}.  
\]
Note that as 
$f$ is homogeneous, $F_f^*$ is the complement of the projective hypersurface
defined by $f=0$.  Let $f = f_1 \cdots f_e$ be a factorization of $f$
into irreducible polynomials.  We note that $(\CC^n\setminus f^{-1}(0))^*$ has 
first homology generated by the meridian circles $\gamma_i$ around $f_i^{-1}(0)\subset \PP^{n-1}$ 
with orientations determined by the complex orientations. %cite Dimca 
%Following D. Cohen and A. Suciu \cite[Corollary 1.4]{article:Cohen}, 
For $\eta^r = 1$ we denote by
$\calV_\eta^f$ the rank one local system on $(\CC\setminus f^{-1}(0))^*= F_f^*$ induced by the homomorphism 
\begin{equation}\label{eqn:Phi}
\Phi_\eta^f : H_1(F_f^*) \to GL_1(\CC) = \CC^*
\end{equation} 
that sends $H_1(\rho)(\gamma_i)$ to $\eta$.  
We define $\calW_\eta$ to be the local system on $\CC^*$ induced by the representation that sends the standard generator 
of $\pi_1(\CC^*)$ to $\eta\in \CC^*$.  
We will also let $H^*(F_f,\CC)_\eta$ denote the
eigenspace of the algebraic PL monodromy of $f$ with eigenvalue
$\eta$.  Lastly, the 
symbol $\#$ will be used as a subscript of a continuous function, and denotes the induced homomorphism defined on  
the fundamental groups.

We will show the following two results:
\begin{lem}\label{lem:fibration}
Let $f$ and $g$ be as in Hypothesis \ref{hyp:assumptions}.  Then there is a fibration $F_{fg}^* \to F_f^* \times F_g^*$ defined by 
\[
[x_1: \ldots :x_n: y_1: \ldots : y_m] \mapsto ([x_1:\ldots : x_n],[y_1:\ldots :y_m])
\]
with fiber $\CC^*$ and trivial geometric monodromy.  
\end{lem}
The Leray spectral sequence associated to this fibration will allow us to prove the following formula.  
\begin{theo}\label{theo:main}
Let $f$ and $g$ satisfy Hypothesis \ref{hyp:assumptions}.  Then,
\[
H^*(F_{fg},\CC)_\eta = H^*(F_f, \CC)_\eta \otimes H^*(F_g, \CC)_\eta
\otimes H^*(\CC^*,\CC).
\] 
\end{theo}

In the statement above, the tensor symbol is used to mean the tensor product 
of vectorspaces graded by cohomological degree.  Namely, if $M_*$ and $N_*$ are graded 
vectorspaces then 
\[
\left( M\otimes N \right)_k = \bigoplus_{i+j=k} M_i \otimes N_j.
\]

It clearly follows from Theorem \ref{theo:main} that the Weil zeta function of the algebraic PL monodromy of $fg$ is
always $1$, recovering N\'emethi's result in our case \cite{article:Nemethi:gen}.

\section{Proof of Theorem \ref{theo:main}}

We will make heavy use of a Theorem of D.Cohen and A. Suciu \cite{article:Cohen}
that we state here.  
\begin{theo}[Cohen, Suciu]\label{theo:Cohen}
Let $f$ be homogeneous and reduced of degree $r$, and pick $\eta \in \CC^*, \eta^r =
1$.  Then
\[
H^*(F_f,\CC)_\eta \cong H^*(F_f^*, \calV^f_\eta).
\]
\qed
\end{theo}

The equation above simply tells us that an eigenspace of the algebraic PL monodromy of $f$ is 
isomorphic to the cohomology of a local system defined on the complement of the projective hypersurface $f = 0$.
This Theorem leads us to consider the cohomology $H^*(F_{fg},
\calV^{fg}_\eta)$.  We establish Lemma \ref{lem:fibration} to aid in the computation of $H^*(F_{fg}, \calV_\eta^{fg})$. 

\begin{lem}
Let $\PP^{n+m-1}$ have coordinates $x_1,\ldots, x_n,y_1,\ldots, y_m$.
Let 
\[
M= \PP^{n+m-1}\setminus (\PP^{n-1} \cup \PP^{m-1})
\] 
be the complement of the projective variety defined by the ideal
\[
(x_1, x_2 , \ldots , x_n) \cdot (y_1,y_2, \ldots , y_n).
\]
Then the map 
\[
\phi:M \to \PP^{n-1}
\times \PP^{m-1}
\]
sending $[x_1 : \ldots : x_n : y_1 : \ldots : y_m]$ to
$([x_1:\ldots : x_n], [y_1 : \ldots : y_m])$ makes $M$ a $\CC^*$-bundle over 
$\PP^{n-1} \times \PP^{m-1}$ with trivial geometric monodromy.
\end{lem} 

\begin{proof}
This map is clearly well-defined.  When we restrict to the chart
defined by $x_j \neq 0$ (resp.\ $y_j \neq 0$) then $\phi$ can be interpreted as the Hopf fibration
applied to the $y$'s (resp.\ $x$'s).  This map is clearly surjective and has trivial geometric 
monodromy as $\PP^{n-1} \times \PP^{m-1}$ is simply connected.
\end{proof}

\begin{proof}[Proof of Lemma \ref{lem:fibration}]
The restriction of $\phi$ to $F_{fg}^*$ has image $F_f^* \times F_g^*$
and has trivial geometric monodromy. 
\end{proof}
When we look at the Leray spectral sequence induced by this
fibration we obtain the following result.  

\begin{theo}\label{theo:ss}
If Hypothesis \ref{hyp:assumptions} holds then there is a spectral sequence 
\begin{equation}\label{eqn:ss}
E_2^{i,j} \Longrightarrow H^{i+j}(F_{fg},\CC)_\eta,
\end{equation}
with $E_2^{i,j} = 0$ for $j\neq 0,1$, and 
\[
E_2^{i,0} \cong E_2^{i,1} \cong 
\bigoplus_{j+k = i} H^j(F_f,\CC)_\eta \otimes H^k(F_g,\CC)_\eta .
\] 
\end{theo}

\begin{proof}
Let $M$ be the complement of $f=0$ in $\CC^n$ and $f: M \to \CC^*$ be the Milnor fibration.  
On page 107 of \cite{article:Cohen} we have the following commutative diagram with exact rows.  
\begin{equation}\label{eqn:diagram}
\xymatrix{
\pi_1(\CC^*) \ar[r]^{\iota}\ar[d]^{\cong} & \pi_1(M) \ar[r] \ar[d]^{f_\#}& \pi_1(M^*) \ar[d]\\
\ZZ \ar[r]^{\times r} & \ZZ \ar[r] & \ZZ}
\end{equation}
The top row of this diagram is part of the homotopy sequence associated to the Hopf fibration restricted to $M$.  
It also follows from \cite{article:Cohen} that if $f=f_1\cdots f_e$ is a factorization of $f$ into distinct 
irreducible polynomials 
and if $a_i$ is the homotopy class of a meridian around $f_i = 0$ with orientation determined by the complex orientations, 
then $f_\#(a_i) = 1$ for all $i= 1,\ldots,e$.  In this way $H_1(M,\ZZ)$ may be identified with the free $\ZZ$ module with 
basis given by the homology classes determined by each of the $a_i$, and $f_*:H_1(M,\ZZ) \to \ZZ$ 
may be identified with the matrix $[1,1,\ldots,1]^T$.  Also note by commutativity that if $\sigma$ is an appropriate choice of a generator of 
$\pi_1(\CC^*)$, then 
we have $f_\# \circ \iota (\sigma) = r$.  Thus in particular $\Phi_\eta^{f} ([\iota(\sigma)]) = \eta^r$, where $[*]$ denotes 
``the homology class determined by''.  

Now recall the fibration from Lemma \ref{lem:fibration}:  
\begin{equation}\label{eqn:phifib}
\xymatrix{
\CC^* \ar[r]^\kappa &  F_{fg}^* \ar[r]^\phi & F_f^* \times F_g^*}.
\end{equation}
Recall further that we consider $F_f^* \times F_g^*$ as a subset of $\PP^{n-1} \times \PP^{m-1}$ with coordinates 
$x_1, \ldots , x_n , y_1 , \ldots , y_m$.  The open subset $y_1 \neq 0$ of $F_f^* \times F_g^*$ can be 
thought of as $M^* \times C$ where $C$ is the 
complement of the hypersurface $g(1,y_2, \ldots , y_m)=0$ in $\AA^{m-1}$.  
In this way, $\phi^{-1}(M^* \times C)$ is $M\times C$ in $\AA^{m+n-1}$.  In fact 
$\phi|_{\phi^{-1}(F_f^* \times C)}$ may be identified with the Hopf fibration applied to the $x$'s, 
\begin{equation}\label{eqn:Cfib}
\xymatrix{
\CC^* \ar[r] & M \times C \ar[r] & M^*\times C}.
\end{equation}
Let $U \subseteq C$ be a contractible subset of $C$ and 
$\psi = \phi |_{\phi^{-1}(F_f^* \times U)}$, and consider the 
restriction of (\ref{eqn:Cfib}), 
\[
\xymatrix{
\CC^* \ar[r]^\lambda & M \times U \ar[r]^\psi & M^* \times U}.
\]
Then by the discussion in the first paragraph of this proof we know that 
\begin{equation}\label{eqn:ttm}
\Phi_\eta^{fg}([\lambda_\# (\sigma ) ]) = \eta^r ,
\end{equation}
where $[*]$ denotes the image of the homology 
class of $*$ under the natural map $H^1(F_f^*\times U ,\ZZ) \to H^1(F_f^* \times F_g^*, \ZZ)$.

Let $\eta^{r+s} = 1$.  We will see that the
spectral sequence of the Theorem is essentially the Leray spectral
sequence 
\begin{equation}\label{eqn:lrss}
H^i(F^*_f \times F_g^*, \RR^j \phi_* (\calV^{fg}_\eta)) \Longrightarrow H^{i+j}(F_{fg}^*,\calV^{fg}_\eta).
\end{equation}
To compute $\RR^\ell\phi_*(\calV^{fg}_\eta)$, we may apply \cite[Proposition 6.4.3]{book:Dimca} to obtain:
\begin{itemize}
\item $\RR^\ell \psi_*(\calV_\eta^{fg}|_{F_f^*\times U})= \pi_1^*(\calV^f_\eta)|_{F_f^*\times U}$ if $\eta^r = 1$ and $\ell = 0,1$.
\item $\RR^\ell \psi_*(\calV_\eta^{fg}|_{F_f^*\times U}) = 0$ otherwise.
\end{itemize}
Here $\pi_1:F_f^* \times F_g^* \to F_f^*$, and $\pi_2: F_f^* \times F_g^* \to F_g^*$ are the natural projections.  
Note that the equation (\ref{eqn:ttm}) calculates what A. Dimca calls the total monodromy operator on 
\cite[p. 210]{book:Dimca}.  A symmetric 
argument holds with $f$ replaced by $g$ (note that $\eta^r = \eta^{-s}$ as $\eta^{r+s} = 1$) and we may conclude:
\begin{itemize}
\item $\RR^\ell \phi_*(\calV^{fg}_\eta) = \pi_1^*(\calV_\eta^f) \otimes 
\pi_2^* (\calV^g_\eta ) := \calV_\eta^f \boxtimes \calV^g_\eta$ if $\eta^r = 1$ and $\ell = 0,1$.
\item $\RR^\ell \phi_*(\calV^{fg}_\eta) = 0$ otherwise.
\end{itemize}

Therefore, if $\eta^r \neq 1, \eta^{r+s} = 1$ the spectral sequence (\ref{eqn:lrss}) is zero.   
Hence $H^*(F_{fg},\CC)_\eta = 0$.  In this case $H^*(F_f, \CC)_\eta$ is also zero and the Theorem is proved.  When $\eta^r = 1$ 
we may now apply the K\"unneth formula \cite[Theorem 4.3.14]{book:Dimca} to obtain that for $j = 0,1$ one has 
\[
H^\ell(F_f^*\times F_g^*, \RR^j \phi_* \calV^{fg}_\eta) 
= \bigoplus_{i+k=\ell} H^i(F_f^*, \calV^f_\eta) \otimes H^k(F_g^* , \calV^f_\eta),
\]
while the left hand side is zero for other $j$.  
These are exactly the $E^{\ell,j}_2$ terms of the spectral sequence of our Theorem, and we know that it 
converges to $H^*(F_{fg}^*, \calV^{fg}_\eta)\cong H^*(F_{fg}^*, \CC)_\eta$.  
This establishes the Theorem for $\eta^{r+s} =1$.

It may be noted that when $\eta^{r+s} \neq 1$ then either $\eta^r$
or $\eta^s$ is not equal to one.  If $\eta^r \neq 1$ then $H^*(F_f,\CC)_\eta$ is zero 
as well as $H^*(F_{fg}, \CC )_\eta$.  If $\eta^s \neq 1$ then $H^*(F_g,\CC)_\eta$ is zero as well as 
$H(F_{fg},\CC)_\eta$ and the 
Theorem follows in these cases.
\end{proof}

We will now concern ourselves with computing $H^i(F_{fg}, \CC)$.  
We first consider a variant of \cite[Theorem 4]{article:Oka} and \cite[Theorem 2]{incollection:Sakamoto}
in our homogeneous case.  We tacitly use the embedding $F_{fg} \subset
\CC^n \times \CC^m$ in the following statement.  

\begin{lem}\label{lem:Sakamoto}
The map $\widehat{f}:F_{fg} \to \CC^*$ defined by $(x,y) \mapsto f(x)$
is a fibration with fiber $F_f \times F_g$.  A geometric PL monodromy of this fibration is, 
\[
(x,y) \mapsto \left(\exp\left(\frac{2\pi i}{r}\right) x, \exp\left(-\frac{2\pi i}{s}\right)y\right).
\]
\end{lem}

\begin{proof}
Since $f$ and 
$g$ are homogeneous we have $f(\exp(\frac{t}{r}) x) = \exp(t)f(x)$ and $g(\exp(\frac{t}{s})y)=\exp(t)g(y)$.  
We also note that if $f(x)g(y)= 1$ then
$f(x) = g(y)^{-1}$.  These two properties prove the Lemma.
\end{proof}

We have a direct consequence, included here for completeness.
\begin{theo}\label{theo:betti}
Let $f$ and $g$ satisfy Hypothesis \ref{hyp:assumptions}, then 
\[
H^\ell(F_{fg},\CC) \cong \bigoplus_{\lambda = \ell -1}^\ell 
\left( \mathop{\bigoplus_{\eta^{r+s}=1}}_{i+j=\lambda} 
H^i(F_f, \CC)_\eta \otimes H^j(F_g,\CC)_{\eta^{-1}} \right)
\]
where the inner sum runs over all possible $\eta$.
\end{theo}

\begin{proof}
The algebraic PL monodromy operator of $\widehat{f}$ is 
\[
T_f \otimes T_g^{-1}: H^*(F_f, \CC) \otimes H^*(F_g, \CC) \to H^*(F_f, \CC) \otimes H^*(F_g, \CC)
\]
where $T_f,T_g$ are the algebraic 
PL monodromy operators of the respective fibers.
Since $T_f,T_g$ are of finite order, they are diagonalizable.  We
let $\{ a_i\}$ (resp.\ $\{ b_j \}$) be a homogeneous basis of $H^*(F_f, \CC)$
(resp.\ $H^*(F_g,\CC)$) that are eigenvectors of  $T_f$ (resp.\ $T_g$),
with eigenvalue $\alpha_i$ (resp.\ $\beta_j$).
In such a case $\{ a_i \otimes b_j \}$ are a basis of eigenvectors for $T_f \otimes
T_g^{-1}$ with eigenvalue $\alpha_i \beta_j^{-1}$.  This shows that 
\[
\RR^\ell \widehat{f}_* (\CC_{F_{fg}}) =
\mathop{\bigoplus_{(\alpha_i,\beta_j)}}_{\deg(a_i) +
  \deg(b_j) = \ell} \calW_{\alpha_i \beta_j^{-1}}.
\]
Ergo, since non-constant rank one local systems on $\CC^*$ have no cohomology
\[
H^p(\CC^*,\RR^\ell\psi_* (\CC_{F_{fg}})) = \mathop{\bigoplus_{\alpha_i \beta_j^{-1}
  = 1}}_{\deg(a_i)+\deg(b_j)=\ell}  \CC
\]
for $p=0,1$, and the left hand side is zero for $p\neq 0,1$.  

Now we may consider the Leray spectral sequence associated with the fibration of Lemma \ref{lem:Sakamoto}.  
Since the base of this fibration is $\CC^*$, the spectral
sequence has only two columns, and thus converges on the
second page.  This yields the Theorem.  
\end{proof}

\begin{proof}[Proof of Theorem \ref{theo:main}]
Theorem \ref{theo:main} now follows from Theorems \ref{theo:ss} and \ref{theo:betti}.  To see this we denote by $\Gamma^\ell_\eta$ the 
vector space $\bigoplus_{i+j=\ell} H^i(F_f,\CC)_\eta \otimes H^j(F_g,\CC)_\eta$.  Theorem \ref{theo:ss} shows the existence of the 
following exact sequence,  
\[
\xymatrix{
\Gamma^{\ell -2}_\eta \ar[r]^{d_2} & \Gamma^\ell_\eta \ar[r] & H^\ell(F_{fg},\CC)_\eta \ar[r] &\Gamma^{\ell-1}_\eta \ar[r]^{d_2} 
& \Gamma^{\ell+1}_\eta}.
\]
By \cite[Proposition 1.1]{article:Cohen}, $H^\ell(F_f,\CC)_\eta \cong H^\ell(F_f,\CC)_{\eta^{-1}}$.  Hence  
\[
H^\ell(F_{fg},\CC) \cong \mathop{\bigoplus_{j = \ell-1}^\ell}_{\eta^{r+s}=1} \Gamma_\eta^j,
\]
by Theorem \ref{theo:betti},
and so the second differential, $d_2$ is zero.
\end{proof}

\section{A Few Examples}

The reader may wish to consult \cite{book:Orlik} for definitions of terms that involve hyperplane arrangements.  
\begin{example}
Let $r = 4, s= 5, f = x_1 x_2 (x_1+x_2)(x_1 + 2 x_2), g = y_1 y_2 (y_1+y_2)(y_1 + 2 y_2)(y_1 + 3y_2)$.  
Note that $f$ and $g$ define generic central line arrangements.  
The Weil zeta function of any generic hyperplane arrangement singularity is presented in 
\cite{article:OrlikRandell}.  To compute the Weil zeta function of any hyperplane arrangement singularity 
one may use \cite[Theorem 9.6]{book:Milnor}, the formula $\chi(F_h) = deg(h)\chi(F_h^*)$, 
and the algorithm \cite[Theorem 5.87(c)]{book:Orlik}.  This method is practical for low 
dimensions and is simple for line arrangements.  Also the Weil zeta 
function of any generic hyperplane arrangement singularity is presented in 
\cite{article:OrlikRandell}.  
Since the 
Milnor fiber of $f$ and $g$ is connected we may easily write down tables expressing the eigenspace decomposition of the 
algebraic (PL) monodromy as follows.  Note that we express $\exp(2\pi i / 5)$ as $\omega$.  

\begin{minipage}{2in}
\begin{center}
$\dim(H^j(F_f,\CC)_\eta)$

\begin{tabular}{|c|c|c|}
\hline
$\eta $\textbackslash$ j$ & $0$  & $1$ \\
\hline
$1$  & $1$ & $3$ \\
$i$ & $0$ & $2$ \\
$-1$ & $0$ & $2$ \\
$-i$ & $0$ & $2$ \\
\hline
\end{tabular}
\end{center}
\end{minipage}
\begin{minipage}{2in}
\begin{center}
$\dim(H^j(F_g,\CC)_\eta)$

\begin{tabular}{|c|c|c|}
\hline
$\eta$\textbackslash$j$ & $0$ & $1$ \\
\hline
$1$ & $1$ & $4$ \\
$\omega$ & $0$ & $3$ \\
$\omega^2$ & $0$ & $3$ \\
$\omega^3$ & $0$ & $3$ \\
$\omega^4$ & $0$ & $3$ \\
\hline
\end{tabular}
\end{center}
\end{minipage}

\noindent Now our Theorem \ref{theo:main} immediately yields the following table for $\dim(H^j(F_{fg},\CC)_\eta)$:
\begin{center}
\begin{tabular}{|c|c|c|c|c|}
\hline
$\eta$\textbackslash$j$ & $0$ & $1$ &$2$&$3$ \\
\hline
$1$ & $1$ & $8$ & $19$ & $12$ \\
\hline
\end{tabular}
\end{center}
where there is a zero for every other $\eta$ and $j$.  In this example even though the algebraic PL monodromy of 
$f$ and $g$ have non-unity eigenvalues the algebraic $PL$ monodromy of $fg$ has one as the only eigenvalue.  
\end{example}
This behavior is not 
uncommon.  Theorem \ref{theo:main} guarantees that $H^*(F_f,\CC)_\eta$ will not contribute to 
$H^*(F_{fg},\CC)$ if $H^*(F_g,\CC)_\eta$ is zero.  We state this observation as the following Corollary of Theorem \ref{theo:main}.  

\begin{carol}
We assume the conditions of Hypothesis \ref{hyp:assumptions}.  
$H^*(F_{fg},\CC)_\eta \neq 0$ if and only if $H^*(F_f,\CC)_\eta \neq 0$ 
and $H^*(F_g,\CC)_\eta \neq 0$.  In particular if $H^*(F_{fg},\CC)_\eta\neq 0$ then $\eta^{\gcd(r,s)} = 1$.  \qed
\end{carol}
Here we give an if and only if condition for the vanishing of $H^j(F_{fg},\CC)_\eta$. 
The second paragraph of remark 3.2 of \cite{incollection:Libgober} only provides the second sentence of this Corollary.  
\begin{example}
Let $r= 3, s=6,f=x_1^3+x_2^3+x_3^3$ and $g = y_1y_2(y_1+y_2)(y_1+2y_2)(y_1+3y_2)(y_1 + 4y_2)$.  The eigenspace decomposition of $f$ is 
discussed in \cite[section 9]{book:Milnor} and $g$ is a line arrangement so we obtain

\begin{minipage}{2in}
\begin{center}
$\dim(H^j(F_f,\CC)_\eta)$

\begin{tabular}{|c|c|c|c|}
\hline
$\eta$\textbackslash$j$& $0$&$1$&$2$\\
\hline 
$1$& $1$ & $0$ & $2$ \\
$\omega^2$&$0$ &$0$&$3$\\
$\omega^4$&$0$ &$0$&$3$\\
\hline
\end{tabular}
\end{center}
\end{minipage}
\begin{minipage}{2in}
\begin{center}
$\dim(H^j(F_g,\CC)_\eta)$

\begin{tabular}{|c|c|c|}
\hline
$\eta$\textbackslash$j$&$0$&$1$ \\
\hline
$1$&$1$&$5$\\
$\omega$&$0$&$4$\\
$\omega^2$&$0$&$4$\\
$\omega^3$&$0$&$4$\\
$\omega^4$&$0$&$4$\\
$\omega^5$&$0$&$4$\\
\hline
\end{tabular}
\end{center}
\end{minipage}

\noindent where $\omega=\exp(2\pi i/ 6)$.  Applying our Theorem \ref{theo:main} yields the following table for 
$\dim(H^j(F_{fg},\CC)_\eta)$.
\begin{center}
\begin{tabular}{|c|c|c|c|c|c|}
\hline
$\eta$\textbackslash$j$&$0$&$1$&$2$&$3$&$4$\\
\hline
$1$&$1$&$6$&$7$&$12$&$10$\\
$\omega^2$&$0$&$0$&$0$&$12$&$12$\\
$\omega^4$&$0$&$0$&$0$&$12$&$12$\\
\hline
\end{tabular}
\end{center}
\end{example}

\section*{Acknowledgments}

The author would like to thank his advisor Uli Walther, and is grateful for helpful conversations with D. Arapura.

\end{document}